\documentstyle[12pt]{article}

\begin{document}

\noindent {\bf REPRESENTATIONS OF THE
$q$-DEFORMED ALGEBRA ${\bf U}_{\bf q}({\bf iso}_{\bf q}({\bf 2}))$}

\vskip 15 pt

{\sl M. Havl\'{i}\v cek}

{\it Department of Mathematics, FNSPE, Czech Technical University}

{\it CZ-120 00, Prague 2, Czech Republic}
\medskip

{\sl A. U. Klimyk}

{\it Institute for Theoretical Physics, Kiev 252143, Ukraine}
\medskip

{\sl S. Po\v sta}

{\it Department of Mathematics, FNSPE, Czech Technical University}

{\it CZ-120 00, Prague 2, Czech Republic}

\vskip 30 pt

\begin{abstract}
An algebra homomorphism $\psi$ from the $q$-deformed
algebra $U_q({\rm iso}_2)$ with generating elements $I$, $T_1$, $T_2$
and defining relations $[I,T_2]_q=T_1$, $[T_1,I]_q=T_2$,
$[T_2,T_1]_q=0$ (where $[A,B]_q=q^{1/2}AB-q^{-1/2}BA$)
to the extension
${\hat U}_q({\rm m}_2)$ of the Hopf algebra $U_q({\rm m}_2)$ is
constructed. The algebra $U_q({\rm iso}_2)$ at $q=1$ leads to the Lie
algebra ${\rm iso}_2 \sim {\rm m}_2$ of the group $ISO(2)$
of motions of the Euclidean plane.
The Hopf algebra $U_q({\rm m}_2)$ is treated as a Hopf $q$-deformation
of the universal enveloping algebra of ${\rm iso}_2$ and is well-known in
the literature.

Not all irreducible representations of
$U_q({\rm m}_2)$ can be extended to representations of the extension
${\hat U}_q({\rm m}_2)$. Composing the homomorphism $\psi$ with irreducible
representations of ${\hat U}_q({\rm m}_2)$ we obtain representations of
$U_q({\rm iso}_2)$. Not all of these representations of $U_q({\rm iso}_2)$
are irreducible. The reducible representations of $U_q({\rm iso}_2)$
are decomposed into irreducible components. In this way we obtain all
irreducible representations of $U_q({\rm iso}_2)$ when $q$ is not a root of
unity. A part of these
representations turns into irreducible representations of the Lie algebra
iso$_2$ when $q\to 1$. Representations of the other part have no
classical analogue.

\end{abstract}

\newpage

\noindent
{\bf I. INTRODUCTION}
\bigskip

Soon after definition of Drinfeld--Jimbo algebras $U_q(g)$,
corresponding to semisimple Lie algebras $g$,
the Hopf algebra $U_q({\rm m}_2)$ was defined [1]
which is treated as a $q$-deformation
of the universal enveloping algebra of the Lie algebra ${\rm iso}_2$ of
the group of motions of the Euclidean plane (the description of this
group, its Lie algebra and their representations see, for example,
in [2], Chap. 4).

However, there is another $q$-deformation of the universal enveloping
algebra $U({\rm iso}_2)$ of the Lie algebra ${\rm iso}_2$ which will
be denoted by $U_q({\rm iso}_2)$. In the general
form (that is, for $U({\rm iso}_n)$) such $q$-deformations were defined in
[3]. The Hopf algebra $U_q({\rm m}_2)$ is related to the well-known
quantum algebra $U_q({\rm sl}_2)$ while the associative algebra
$U_q({\rm iso}_2)$ is connected with the nonstandard $q$-deformation
$U_q({\rm so}_3)$ of the universal enveloping algebra $U({\rm so}_3)$
which is sometimes called the Fairlie algebra.

It is known that the theory of representations of the associative
algebra $U_q({\rm so}_3)$ is richer than that of the algebra
$U_q({\rm sl}_2)$ [4--6]. It was shown recently [7] that the theory of
representations of the algebra $U_q({\rm iso}_2)$ is also richer
than that of the algebra $U_q({\rm m}_2)$. In particular, the algebras
$U_q({\rm so}_3)$ and $U_q({\rm iso}_2)$ have irreducible representations
of nonclassical type (that is, representations which have no limit at
$q\to 1$). The paper [7] is devoted to study of irreducible
$*$-representations of the algebra $U_q({\rm iso}_2)$ equipped with
$*$-structures. Irreducible representations of $U_q({\rm iso}_2)$
of the classical type are given in [8].

The aim of the present paper is to study irreducible
representations of $U_q({\rm iso}_2)$ when this algebra is not
equipped with some $*$-structure and to clarify why irreducible
representations of $U_q({\rm iso}_2)$ of
the nonclassical type appear. We do this in the same
way as in the case of representations of the algebra $U_q({\rm so}_3)$
in [6]. Namely, we relate the algebra $U_q({\rm iso}_2)$ with the
extension ${\hat U}_q({\rm m}_2)$ of the Hopf algebra
$U_q({\rm m}_2)$. This allows us to obtain representations of
$U_q({\rm iso}_2)$ from those of the extended
algebra ${\hat U}_q({\rm m}_2)$. We prove that if  $q$ is not a root
of unity, then irreducible
representations obtained in this way exhaust, up to equivalence, all
irreducible representations of $U_q({\rm iso}_2)$.
\bigskip

\noindent
{\bf II. THE ALGEBRAS $U_q({\rm iso}_2)$ AND ${\hat U}_q({\rm m}_2)$}
\medskip

The algebra $U_q({\rm iso}_2)$ is obtained by a $q$-deformation of the
standard commutation relations
$$
[I,T_2]=T_1,\ \ \ [T_1,I]=T_2,\ \ \ [T_2,T_1]=0
$$
of the Lie algebra iso$_2$. So, $U_q({\rm iso}_2)$ is defined [7, 8] as the
complex associative algebra with unit element generated by the elements
$I$, $T_1$, $T_2$ satisfying the defining relations
$$
[I,T_2]_q:= q^{1/2}IT_2-q^{-1/2}T_2I=T_1,  \eqno (1)$$
$$
[T_1,I]_q:= q^{1/2}T_1I-q^{-1/2}IT_1=T_2, \eqno (2) $$
$$
[T_2,T_1]_q :=q^{1/2}T_2T_1-q^{-1/2}T_1T_2=0. \eqno (3)
$$
Note that the elements $T_2$ and $T_1$ of the algebra $U_q({\rm iso}_2)$
do not commute (as it is a case in the algebra ${\rm iso}_2$; these
elements correspond to shifts along the axes of the plane). We say that
they $q$-commute, that is, $q^{1/2}T_2T_1-q^{-1/2}T_1T_2=0$.
This means that they generate the associative algebra determining
the quantum plane.

Unfortunately, a Hopf algebra structure is not known on $U_q({\rm iso}_2)$.
However, it can be embedded into the Hopf algebra $U_q({\rm isl}_2)$ as a
Hopf coideal. (The algebra $U_q({\rm isl}_2)$ is the $q$-deformation
of the universal enveloping algebra $U({\rm isl}_2)$ of the Lie
algebra ${\rm isl}_2$ of the inhomogeneous Lie group $ISL(2)$).

The relations (1)--(3) lead to the Poincar\'e--Birkhoff--Witt theorem
for the algebra $U_q({\rm iso}_2)$. This theorem can be
formulated as:
\medskip

\noindent {\bf Proposition 1.}
{\it The elements $T_1^jT_2^kI^l$, $j,k,l=0,1,2,\cdots $,
form a basis of the linear space $U_q({\rm iso}_2)$.}
\medskip

Indeed, by using
the relations (1)--(3) any product of the elements $I$, $T_2$, $T_1$
can be reduced to a sum of the elements $T_1^jT_2^kI^l$ with complex
coefficients. Using the diamond lemma [9] (or its special case from Subsect.
4.1.5 in [10]) it is proved that these elements are linear independent.
This proves Proposition 1.
\medskip

Note that by (1) the element $T_1$ is not independent: it is determined
by the elements $I$ and $T_2$. Thus, the algebra $U_q({\rm iso}_2)$
is generated by $I$ and $T_2$, but now instead of quadratic relations
(1)--(3) we must take the relations
$$
I^2T_2-(q+q^{-1})IT_2I+T_2I^2=-T_2, \eqno (4)$$
$$
IT_2^2-(q+q^{-1})T_2IT_2+T_2^2I=0, \eqno (5)
$$
which are obtained if we substitute the expression (1) for $T_1$ into
(2) and (3). The equation $q^{1/2}IT_2-q^{-1/2}T_2I=T_1$ and the
relations (4) and (5) restore the relations (1)--(3).

Note that the relation (5) is a relation of Serre's type
in the definition of quantum algebras by V. Drinfeld and M. Jimbo.
The relation (4) differs from Serre's relation by
appearance of non-vanishing right hand side.

It is known that the element $C=T_1^2+T_2^2$ from the universal
enveloping algebra $U({\rm iso}_2)$
belongs to the center of this algebra. The analogue of this element in
$U_q({\rm iso}_2)$ is the element
$C_q=\frac 12 (T_1T'_1+T'_1T_1)+{\frac 12} (q+q^{-1}) T_2^2,$
where $T'_1=q^{-1/2}IT_2-q^{1/2}T_2I$ (see [8]), that is
$[C_q, X]:=C_qX-XC_q=0$ for all $X\in U({\rm iso}_2)$.
This element can be reduced according to Proposition 1 to the form
$$
C_q=q^{-1}T_1^2+qT_2^2+q^{-3/2}(1-q^2)T_1T_2I. \eqno (6)
$$

The algebra $U_q({\rm iso}_2)$ is closely related to (but not coincides
with) the quantum algebra $U_q({\rm m}_2)$. The last algebra is generated
by the elements $q^H$, $q^{-H}$, $E$, $F$ satisfying the relations
$$
q^Hq^{-H}=q^{-H}q^H=1, \ \ q^HEq^{-H}=qE, \ \ q^HFq^{-H}=q^{-1}F, \ \
[E,F]:= EF-FE=0. \eqno (7)
$$

In order to relate the algebras $U_q({\rm iso}_2)$ and $U_q({\rm m}_2)$
we need to extend $U_q({\rm m}_2)$ by the elements
$(q^kq^H+q^{-k}q^{-H})^{-1}$, $k\in {\bf Z}$, in the sense of [11].
This extension ${\hat U}_q({\rm m}_2)$ is defined as the associative
algebra (with unit element) generated by the elements
$$
q^H,\ \  \ q^{-H},\ \ \ E,\ \ \ F,\ \ \ (q^kq^H+q^{-k}q^{-H})^{-1},
\ \ k\in {\bf Z},
$$
satisfying the defining relations (7)
of the algebra $U_q({\rm m}_2)$ and
the following natural relations:
$$
(q^kq^H+q^{-k}q^{-H})^{-1}(q^kq^H+q^{-k}q^{-H})=
(q^kq^H+q^{-k}q^{-H})(q^kq^H+q^{-k}q^{-H})^{-1}=1, \eqno (8) $$
$$
q^{\pm H}(q^kq^H+q^{-k}q^{-H})^{-1}=(q^kq^H+q^{-k}q^{-H})^{-1}q^{\pm H} ,
\eqno (9) $$
$$
(q^kq^H+q^{-k}q^{-H})^{-1}E=E(q^{k+1}q^H+q^{-k-1}q^{-H})^{-1}, \eqno (10) $$
$$
(q^kq^H+q^{-k}q^{-H})^{-1}F=F(q^{k-1}q^H+q^{-k+1}q^{-H})^{-1}. \eqno (11)
$$
\medskip

\noindent
{\bf III. THE ALGEBRA HOMOMORPHISM
$U_q({\rm iso}_2)\to {\hat U}_q({\rm m}_2)$}
\medskip

The aim of this section is to give (in an explicit form) the homomorphism
of the algebra $U_q({\rm iso}_2)$ to ${\hat U}_q({\rm m}_2)$.
\medskip

\noindent
{\bf Proposition 2.} {\it There exists a unique algebra homomorphism
$\psi : U_q({\rm iso}_2)\to {\hat U}_q({\rm m}_2)$ such that
$$
\psi (I_1)=\frac {\rm i}{q-q^{-1}} (q^H-q^{-H}), \eqno (12) $$
$$
\psi (I_2)=(E-F) (q^H+q^{-H})^{-1}, \eqno (13) $$
$$
\psi (I_3)=({\rm i} q^{H-1/2}E+{\rm i}q^{-H-1/2}F) (q^H+q^{-H})^{-1},
\eqno (14)
$$
where $q^{H+a}:= q^Hq^a$ for $a\in {\bf C}$.}
\medskip

{\sl Proof.} In order to prove this proposition we have to show that
the defining relations
$$
q^{1/2}\psi (I)\psi (T_2)-q^{-1/2}\psi (T_2)\psi (I)=\psi (T_1), $$
$$
q^{1/2}\psi (T_1)\psi (I)-q^{-1/2}\psi (I)\psi (T_1)=\psi (T_2),
$$
$$
q^{1/2}\psi (T_2)\psi (T_1)-q^{-1/2}\psi (T_1)\psi (T_2)=0.
\eqno (15)
$$
of $U_q({\rm iso}_2)$ are satisfied.
Let us prove the relation (15). (Other relations are proved similarly.)
Substituting the expressions (12)--(14) for $\psi (I)$, $\psi (T_2)$,
$\psi (T_1)$ into (15) we obtain (after multiplying both
sides of equality by $(q^H+q^{-H})$ on the right) the relation
$$
q(E-F)Eq^H(qq^H+q^{-1}q^{-H})^{-1}+q(E-F)Fq^{-H}(q^{-1}q^H+qq^{-H})^{-1}- $$
$$
-qE^2q^H(qq^H+q^{-1}q^{-H})^{-1}-q^{-1}FEq^{-H}(qq^H+q^{-1}q^{-H})^{-1}+ $$
$$
+q^{-1}EFq^H(q^{-1}q^H+qq^{-H})^{-1}+qF^2q^{-H}(q^{-1}q^H+qq^{-H})^{-1}
=0.
$$
The formula (15) is true if and only if this relation is correct. We
multiply both its sides by $(qq^H+q^{-1}q^{-H})(q^{-1}q^H+qq^{-H})$
on the right and obtain the relation in the algebra $U_q({\rm m}_2)$
(that is, without the expressions $(q^kq^H+q^{-k}q^{-H})^{-1}$). This
relation is easily verified by using the defining relations (7)
of the algebra $U_q({\rm m}_2)$. Proposition is proved.
\bigskip

\noindent
{\bf IV. DEFINITION OF REPRESENTATIONS OF $U_q({\rm m}_2)$ AND
$U_q({\rm iso}_2)$}
\medskip

From this point we assume that $q$ is not a root of unity.
Let us define representations of the algebras
$U_q({\rm m}_2)$ and $U_q({\rm iso}_2)$.
\medskip

\noindent {\bf Definition.} {\it
By a representation $\pi$ of $U_q({\rm m}_2)$ (respectively
$U_q({\rm iso}_2)$) we mean a homomorphism of
$U_q({\rm sl}_2)$ (respectively $U_q({\rm iso}_2)$)
into the algebra of linear operators (bounded or unbounded) on a Hilbert
space ${\cal H}$, defined on an everywhere dense invariant subspace
${\cal D}$, such that the operator $\pi (q^H)$ (respectively the
operator $\pi (I)$) can be diagonalized, has a discrete spectrum
(with finite multiplicities of spectral points if $\pi$ is irreducible)
and its eigenvectors belong to ${\cal D}$.
Two representations $\pi$ and $\pi '$ of $U_q({\rm m}_2)$ (of
$U_q({\rm iso}_2)$) on spaces
${\cal H}$ and ${\cal H}'$, respectively, are called (algebraically)
equivalent if there exist everywhere dense invariant subspaces
$V\subset {\cal H}$ and $V'\subset {\cal H}'$ and a one-to-one linear
operator $A: V\to V'$ such that $A\pi (a) v=\pi '(a)A v$ for all
$a\in U_q({\rm m}_2)$ (respectively, for all $a\in U_q({\rm iso}_2)$)
and $ v\in V$.}
\medskip

\noindent {\sl Remark.} Note that the element $I\in U_q({\rm iso}_2)$
corresponds to the homogeneous part of the motion group $ISO(2)$.
As in the classical case, it is natural to demand in the definition of
representations of $U_q({\rm iso}_2)$ that the operator $\pi (I)$ has
a discrete spectrum (with finite multiplicities of spectral points for
irreducible representations $\pi$). Such representations correspond to
Harish-Chandra modules of Lie algebras.  Note that irreducible
$*$-representations of $U_q({\rm iso}_2)$  without a requirement that
$\pi (I)$ has a discrete spectrum were studied in [7]. It was shown there
that the classification of irreducible $*$-representations by
self-adjoint operators in this case is equivalent to the classification of
arbitrary families of bounded self-adjoint operators. The classification of
irreducible representations (not obligatory $*$-representations) in this
case turn into unsolved problem.
\medskip

The algebra $U_q({\rm m}_2)$ has the following non-trivial
irreducible representations:
\medskip

(a) one-dimensional representations $\pi _\sigma$, $\sigma \in {\bf C}$,
$\sigma \ne 0$, determined by the formulas $\pi _\sigma (q^H)=\sigma$,
$\pi _\sigma (E)=\pi _\sigma (F)=0$;

(b) infinite dimensional representations $\pi _{rs}$, $r,s\in {\bf C}$,
$r,s\ne 0$, acting on the Hilbert space ${\cal H}$ with a basis
$| m\rangle$, $m \in {\bf Z}$, by the formulas
$$
\pi _{rs} (q^H)|m\rangle =sq^m |m\rangle ,\ \ \
\pi _{rs} (E)| m\rangle =r |m+1 \rangle ,\ \ \
\pi _{rs} (F)| m\rangle =r |m-1 \rangle ,\ \ \ m\in {\bf Z} .  \eqno (16)
$$
We take ${\cal D}={\rm lin}\, \{ |m\rangle \, |\, m\in {\bf Z}\} $.
A direct verification shows:
\medskip

\noindent {\bf Proposition 3.} {\it
The representations $\pi _{rs}$ and $\pi _{r's'}$ ($r,s,r',s'\in {\bf C}
\backslash \{ 0\}$) are equivalent if and only if
$r=\pm r'$ and $s'=q^ns$ for some $n\in {\bf Z}$.}
\medskip

Repeating the reasonings of Sect. 5.2 from [10] we easily prove
\medskip

\noindent {\bf Proposition 4.} {\it
Every irreducible representation of $U_q({\rm m}_2)$ is equivalent
to one of the representations (16) or is one-dimensional.}
\medskip

Note that for $q\to 1$ the representations $\pi _{rs}$ of
$U_q({\rm m}_2)$ turn into irreducible representations of the universal
enveloping algebra $U({\rm m}_2)$, that is, all irreducible
representations of $U_q({\rm m}_2)$ are deformations of the corresponding
irreducible representations of $U({\rm m}_2)$.

We try to extend representations $\pi _{r,s}$ of $U_q({\rm m}_2)$ to
representations of the extension
${\hat U}_q({\rm m}_2)$ by using the relation
$$
\pi ((q^kq^H+q^{-k}q^{-H})^{-1}):=(q^k\pi (q^H)+q^{-k}\pi (q^{-H}))^{-1},
\ \ \ k\in {\bf Z}. \eqno (17)
$$
Clearly, only those irreducible representations $\pi _{r,s}$ of
$U_q({\rm m}_2)$
can be extended to ${\hat U}_q({\rm m}_2)$ for which the operators
$q^k\pi (q^H)+q^{-k}\pi (q^{-H})$ are invertible.
From formulas (16) it is
clear that these operators are always invertible for the representations
$\pi _{rs}$, $s\ne \pm {\rm i}q^n$, $n\in {\bf Z}$.
(For the representations
$\pi _{rs}$, $s= \pm {\rm i}q^n$ for some $n\in {\bf Z}$,
some of these operators are not invertible since they have zero
eigenvalue.) Denoting the extended representations by the same symbols
$\pi _{rs}$, we can formulate the following statement:
\medskip

\noindent
{\bf Proposition 5.} {\it The algebra ${\hat U}_q({\rm m}_2)$
has the infinite dimensional representations $\pi _{rs}$,
$r,s\in {\bf C}\backslash \{ 0\}$, $s\ne \pm {\rm i}q^n$ for all
$n\in {\bf Z}$, given by the relations (16) and (17). The representations
$\pi _{rs}$ and $\pi _{r',s'}$ ($r,s,r',s'\in {\bf C}\backslash \{ 0\}$,
$s,s' \ne \pm {\rm i}q^n$ for all $n\in {\bf Z}$)
are equivalent if and only if $r=\pm r'$ and $s'=q^ms$ for some
$m\in {\bf Z}$.
Any irreducible representation of ${\hat U}_q({\rm sl}_2)$ is equivalent to
the representation $\pi _{r,s}$ for some $r,s$ or is a one-dimensional
represenration.}
\bigskip

\noindent
{\bf V. IRREDUCIBLE REPRESENTATIONS OF $U_q({\rm iso}_2)$}
\medskip

If $\pi$ is a representation of the algebra ${\hat U}_q({\rm m}_2)$ on a
space ${\cal H}$, then the mapping $R: U_q({\rm iso}_2)\to {\cal H}$
defined as the composition $R=\pi\circ \psi$, where $\psi$ is the
homomorphism from Proposition 1, is a (not necessary irreducible)
representation of $U_q({\rm iso}_2)$.

Let us consider the representations
$$
R_{rs}=\pi_{rs}\circ \psi
$$
of $U_q({\rm iso}_2)$, where $\pi _{rs}$
are the irreducible representations
of ${\hat U}_q({\rm m}_2)$ from Proposition 5.
Using formulas (16) and (12)--(14) we find that
$$
R_{rs}(I)| m\rangle ={\rm i}\frac {sq^m-s^{-1}q^{-m}}{q-q^{-1}} |m\rangle ,
\eqno (18)$$
$$
R_{rs}(T_2)| m\rangle =\frac r{sq^m+s^{-1}q^{-m}} \{  |m+1\rangle
+ |m-1\rangle \} ,  \eqno (19) $$
$$
R_{rs}(T_1)| m\rangle =\frac {{\rm i}q^{1/2}r}{sq^m+s^{-1}q^{-m}}
\{ sq^m |m+1\rangle -q^{-m}s^{-1} |m-1\rangle \} . \eqno (20)
$$
We consider that these operators are defined on the invariant subspace
${\cal D}\subset {\cal H}$ which is the span of the basis vectors
$|m\rangle$. Thus we proved the following
\medskip

\noindent
{\bf Proposition 6.} {\it Let $r,s\in {\bf C}\backslash \{ 0\}$,
$s\ne \pm {\rm i}q^{n}$ for all $n\in {\bf Z}$. Then the formulas
(18)--(20) form a representation $R_{rs}$ of the algebra
$U_q({\rm iso}_2)$.}
\medskip

We also have
\medskip

\noindent
{\bf Proposition 7.} {\it The representations $R_{rs}$ of Proposition 6
are irreducible if $s\ne \pm {\rm i}q^{m+1/2}$, $m\in {\bf Z}$.}
\medskip

{\sl Proof.} Let $\{ a\}:=(sq^a-a^{-1}q^{-a})/(q-q^{-1})$.
To prove this proposition we first note that since $q$ is not a root of
unity and $s\ne \pm {\rm i}q^{m+1/2}$, $m\in {\bf Z}$, the eigenvalues
${\rm i}\{ m\} $,
$m=0,\pm 1, \pm 2,\cdots$, of the operator $R_{rs}(I)$
are pairwise different.

Let $V\subset {\cal D}$ be an invariant subspace of
the representation $R_{rs}$. We need to show that $V={\cal D}$.
Let $ v=
\sum _{m_i} \alpha _i| m_i\rangle \in V$, where $| m_i\rangle $ are
eigenvectors of $R_{rs}(I)$ which are basis vectors of ${\cal H}$.
(Note that the sum is finite since $ v\in
{\cal D}$.) Let us prove that $|m_i\rangle \in V$. We prove
this for the case when ${\bf v}=\alpha _1|m_1\rangle +\alpha _2|m_2\rangle$.
(The case of more number of summands is proved similarly.) We have
$$
v':= R_{rs}(I) v={\rm i}\alpha _1\{ m_1\} |m_1\rangle +{\rm i}\alpha _2
\{ m_2\} |m_2\rangle .
$$.
Since $v,v'\in V$, one derive that
$$
{\rm i}\{ m_1\} v- v'={\rm i}\alpha _2(\{ m_1\} -
\{ m_2\} )| m_2\rangle \in V.
$$
Since $\{ m_1\} \ne \{ m_2\}$, then  $|m_2\rangle \in V$ and hence
$|m_1\rangle \in V$.

In order to prove that $V={\cal D}$ we obtain from (18) and (19) that
$$
\{ R_{rs}(T_1)-{\rm i} sq^{m_2+1/2}R_{rs}(T_2)\} |m_2\rangle
={\rm i}rq^{1/2}|m_2-1\rangle ,   $$
$$
\{ R_{rs}(T_1)+{\rm i}s^{-1}q^{-m_2+1/2}R_{rs}(T_2)\} |m_2\rangle
={\rm i}rq^{1/2}|m_2+1\rangle .   (20)
$$
It follows from these relations that $V$ contains the vectors
$|m_2-1\rangle , |m_2-2\rangle ,\cdots $ and the vectors
$|m_2+1\rangle ,|m_2+2\rangle ,\cdots $. This means that
$V={\cal D}$ and the representation $R_{rs}$ is irreducible.
Proposition is proved.
\medskip

Note that the representations $R_{rs}$ of Proposition 7 turn
into irreducible representations of the universal enveloping algebra
$U({\rm iso}_2)$ when $q\to 1$. For this reason, they are called
representations of the classical type.

Using Proposition 5 it is easy to show that the representations
$R_{rs}$ and $R_{r's'}$ of Proposition 7 are equivalent if and only if
$r'=\pm r$ and $s =q^ms$ for some $ m\in {\bf Z}$.
\medskip

\noindent
{\bf Proposition 8.} {\it
Let $r\in {\bf C}\backslash \{ 0\}$ and
$s=\varepsilon {\rm i}q^{m+1/2}$, where $m\in {\bf Z}$ and
$\varepsilon \in \{1,-1\}$. Then
the representation $R_{rs}$ is reducible.}
\medskip

{\it Proof.} The eigenvalues of the operator $R_{rs}(I)$ are
$$
-\varepsilon \frac{q^n+q^{-n}}{q-q^{-1}},\ \ \ n=\pm \frac{1}{2},
\pm \frac{3}{2},\pm \frac{3}{2}, \cdots,
$$
that is, every spectral point has
multiplicity 2. The pairs of vectors $|-m+j\rangle$ and
$|-m-j-1\rangle$, $j=0,1,2,\cdots$,
are of the same eigenvalue. Let us
define two subspaces $V_1$ and $V_{-1}$ by the formulas
$V_{\tilde\varepsilon}:=
{\rm lin}\, \{|j\rangle_{\tilde\varepsilon}\, |\, j=0,1,2,\cdots \}$,
where
$$
|j\rangle_{\tilde\varepsilon}:=|-m+j\rangle+
\tilde\varepsilon {\rm i}(-1)^{j}|-m-j-1\rangle,\ \ \
j=0,1,2,\cdots . \eqno (21)
$$
A direct calculation show that for ${\tilde\varepsilon}=1$ and for
${\tilde\varepsilon}=-1$ we have
$$
R_{rs}(I)
|j\rangle_{\tilde\varepsilon}=
-\varepsilon \frac{q^{j+1/2}+q^{-j-1/2}}{q-q^{-1}}
|j\rangle_{\tilde\varepsilon},\ \ \ j=0,1,2,\cdots , \eqno (22)
$$
$$
R_{rs}(T_2)|0\rangle_{\tilde\varepsilon}
=-\frac{r}{q^{1/2}-q^{-1/2}}
\Bigl(\tilde\varepsilon |0\rangle_{\tilde\varepsilon}
+{\rm i} |1\rangle_{\tilde\varepsilon}\Bigr), \eqno (23) $$
$$
R_{rs}(T_2)|j\rangle_{\tilde\varepsilon}
=-\varepsilon\frac{ {\rm i} r}{q^{j+1/2}-q^{-j-1/2}}
\Bigl(|j+1\rangle_{\tilde\varepsilon}+
|j-1\rangle_{\tilde\varepsilon}\Bigr),\ \ \ j=1,2,3,\cdots , \eqno (24) $$
$$
R_{rs}(T_1)|0\rangle_{\tilde\varepsilon}
=\frac{r}{q^{1/2}-q^{-1/2}}
\Bigl(\tilde\varepsilon|0\rangle_{\tilde\varepsilon}
+{\rm i} q|1\rangle_{\tilde\varepsilon}\Bigr), \eqno (25) $$
$$
R_{rs}(T_1)|j\rangle_{\tilde\varepsilon}
=\frac{ {\rm i} r}{q^{j+1/2}-q^{-j-1/2}}
\Bigl(q^{j+1}|j+1\rangle_{\tilde\varepsilon}+
q^{-j}|j-1\rangle_{\tilde\varepsilon}\Bigr),\ \ \ j=1,2,3,\cdots .\eqno (26)
$$
These formulas show that the subspaces $V_1$ and $V_{-1}$ are
invariant with respect to the representation $R_{rs}$, that is
this representation is reducible. Proposition is proved.
\medskip

Let us denote the restrictions of the representation $R_{rs}$
from Proposition 7 to the invariant subspases $V_1$ and $V_{-1}$ by
$R_r^{\varepsilon,1}$ and $R_r^{\varepsilon,-1}$, respectively.
It is seen from the formulas (22)--(26) that
the operators are independent of $s$ and the index $s$ is
ommitted. These formulas show that $R_{rs}$,
$s=\varepsilon {\rm i}q^{m+1/2}$, is the direct sum
of the representations $R_r^{\varepsilon,1}$ and
$R_r^{\varepsilon,-1}$.
\medskip

\noindent
{\bf Proposition 9.} {\it The representations
$R_r^{\varepsilon,\tilde\varepsilon}$ and
$R_{r'}^{\varepsilon',\tilde\varepsilon'}$ are equivalent if
$(r,\varepsilon,\tilde\varepsilon)=
(-r',\varepsilon',-\tilde\varepsilon')$.}
\medskip

Proposition is easily proved by using Proposition 5 and formulas
(22)--(26).
\medskip

\noindent
{\bf Theorem 1.} (a) {\it Let $r$ and $r'$ be nonzero complex numbers
such that ${\rm Re} r>0$ and ${\rm Re} r'>0$, and let
$\varepsilon,\tilde\varepsilon,\varepsilon',\tilde\varepsilon'
\in \{1,-1\}$. If $(\varepsilon,\tilde\varepsilon,r)\ne
(\varepsilon',\tilde\varepsilon',r')$, then
the representations $R_r^{\varepsilon,\tilde\varepsilon}$ and
$R_{r'}^{\varepsilon',\tilde\varepsilon'}$
are irreducible and nonequivalent.}

(b) {\it Let $r \in {\bf C}\backslash \{0\}$, $\varepsilon,\tilde\varepsilon
\in \{1,-1\}$, and let $r',s'\in {\bf C}\backslash \{0\}$,
$s'\ne \pm {\rm i} q^{m+1/2}$ for all $m \in {\bf Z}$. Then
the representations $R_r^{\varepsilon,\tilde\varepsilon}$ and
$R_{r's'}$ are nonequivalent.}
\medskip

{\sl Proof.} The irreducibility is proved in the same way
as in Proposition 7. In order to prove a nonequivalence we note
that the spectrum of the operator $R(I)$ for any of
the representations $R^{+,+}_r$, $R^{-,+}_r$, ${\rm Re}\, r> 0$,
does not coincide with that of any of the representations $R^{+,-}_r$,
$R^{-,-}_r$, ${\rm Re}\, r> 0$. Therefore, any of the representations
$R^{+,+}_r$, $R^{-,+}_r$, ${\rm Re}\, r> 0$, cannot be equivalent
to some of the representations $R^{+,-}_r$, $R^{-,-}_r$,
${\rm Re}\, r> 0$.

The operators $R^{\epsilon _1,\epsilon _2}_r(T_2)$, $\epsilon _1,
\epsilon _2=+,-$, are trace class operators. Their traces are nonzero
(there exists
only one nonzero diagonal matrix element with respect to the basis
$\{ |m\rangle '\}$ or the basis $\{ |m\rangle ''\}$).
Since for ${\rm Re}\, r>0$ and ${\rm Re}\, r'>0$, $r\ne r'$,
we have ${\rm Tr}\, R^{+,+}_r(T_2)\ne {\rm Tr}\, R^{-,+}_{r'}(T_2)$, then
any of the representations $R^{+,+}_r$, ${\rm Re}\, r>0$, cannot be
equivalent to some of the representations $R^{-,+}_{r}$, ${\rm Re}\, r>0$.
It is proved similarly that
any of the representations $R^{+,-}_r$, ${\rm Re}\, r>0$, cannot be
equivalent to some of the representations $R^{-,-}_{r'}$, ${\rm Re}\, r'>0$.
This prove the assertion (a). The assertion (b) is proved similarly.
Proposition is proved.
\medskip

Representations of Theorem 1 have no classical limit since at $q\to 1$ the
denominators in (22)--(26) turn into zero. For this reason, these
representations are called representations of non-classical type. There
are no analogues of such representations for the Lie algebra ${\rm iso}_2$.

\medskip

\noindent
{\bf Theorem 2.} {\it Every irreducible representation of
$U_q({\rm iso}_2)$ is equivalent to one of the representations
$R_{rs}$ and $R_{r}^{\varepsilon,\tilde\varepsilon}$
or is one-dimensional.
This means that the representations $R_{rs}$ and
$R_{r'}^{\varepsilon,\tilde\varepsilon}$, $r,s,r' \in {\bf C}\backslash
\{0\}$, $\varepsilon,\tilde\varepsilon \in \{1,-1\}$ defined by the
relations (18)--(20) and (22)--(26), respectively,
exhaust (up to equivalence) all irreducible representations
of $U_q({\rm iso}_2)$.}
\medskip

{\it Proof.} Let $R$ be an irreducible representation
of $U_q({\rm iso}_2)$.
Then it follows from the definition of representations
that $R(I)$ has some eigenvector
$|0\rangle$. Thus there exists $s \in {\bf C}$, $s\ne 0$, such that
$$
R(I)|0\rangle={\rm i} [0]_{q,s}|0\rangle,
$$
where $[m]_{q,s}:=(s q^m-s^{-1} q^{-m})/(q-q^{-1})$.
Since $R$ is irreducible there exists a complex number $C$ such that
$R(C_q)=C$ (see (6)). We define recursively the vectors
$$
|j+1\rangle:=R({\rm i} T_1-s^{-1} q^{-j+1/2} T_2) |j\rangle,
\ \ \ j= 0,1,2,\cdots , \eqno (27) $$
$$
|j-1\rangle:=R({\rm i} T_1+s q^{j+1/2} T_2) |j\rangle,
\ \ \ j=0,-1,-2,\cdots . \eqno (28)
$$
Some of these vectors may be linear dependent or be equal to 0.
It follows from (1)--(3) and (6) that
$$
R(I) |j\rangle=i [j]_{q,s} |j\rangle, \ \ \ j \in {\bf Z}, \eqno (29) $$
$$
R({\rm i} T_1 +s q^{j+3/2} T_2) |j+1\rangle=-Cq|j\rangle ,
\ \ \ j=0,1,2,\cdots , \eqno (30) $$
$$
R({\rm i} T_1 -s^{-1} q^{-j+3/2} T_2) |j-1\rangle=-Cq|j\rangle,
\ \ \ j= 0,-1,-2,\cdots . \eqno (31)
$$
As a sample, we prove the relation (30) for $j\ge 0$:
$$
R({\rm i} T_1 +s q^{j+3/2} T_2) |j+1\rangle=
R({\rm i} T_1 +s q^{j+3/2} T_2)
R({\rm i} T_1 -s^{-1} q^{-j+1/2} T_2) |j\rangle=  $$
$$
=R(- T^2_1 +{\rm i}s q^{j+3/2} T_2T_1
-{\rm i}s^{-1} q^{-j+1/2}T_1 T_2-q^2T_2^2) |j\rangle=  $$
$$
=qR(- q^{-1}T^2_1 -qT^2_2 +{\rm i}s q^{j-3/2} T_1T_2
-{\rm i}s^{-1} q^{-j-1/2}T_1 T_2) |j\rangle=  $$
$$
=qR(- q^{-1}T^2_1 -qT^2_2 -q^{-3/2}(1-q^2)T_1T_2I)|j\rangle
=-qR(C_q)|j\rangle =-qC|j\rangle.
$$

We obtain from (27) and (30) that
$$
R(T_2)|j\rangle =- (s^{-1} q^{-j+1/2}+sq^{j+1/2})^{-1}
(|j+1\rangle +Cq|j-1\rangle ) , \eqno (32) $$
$$
{\rm i} R(T_1)|j\rangle =\frac{sq^{j+1/2}}{ s^{-1}q^{-j+1/2}+sq^{j+1/2}}
|j+1\rangle +Cq\left( \frac{sq^{j+1/2}}{ s^{-1}q^{-j+1/2}+sq^{j+1/2}}
\right) |j-1\rangle . \eqno (33)
$$
Let us now consider two cases: (a) $C=0$, (b) $C\ne 0$.

(a) $C=0$. The formulas (32) and (33) in this case give
$$
R(T_2)|j\rangle =- (s^{-1} q^{-j+1/2}+sq^{j+1/2})^{-1}|j+1\rangle ,
 \eqno (34)   $$
$$
{\rm i} R(T_1)|j\rangle =\frac{sq^{j+1/2}}{ s^{-1}q^{-j+1/2}+sq^{j+1/2}}
|j+1\rangle  . \eqno (35)
$$
If the set $\{|j\rangle \, |\, j=0,1,2,\cdots \}$ is
linear independent, it follows from (34) and (35) that
${\rm lin}\, \{|j\rangle,|j+1\rangle,\cdots \}$ is
invariant subspace for any $j=1,2,3,\cdots $. Thus the representation
is either reducible or one-dimensional.

Now let there exist $l \in {\bf N}$ such that $|l\rangle$ is linear
dependent on linear independent vectors $|0\rangle,|1\rangle, \cdots ,
|l-1\rangle$. Since the sequence of numbers
$[j]_{q,s}, j \in {\bf Z}$, does
not contain 3 equal elements, the only possible case is
$|l\rangle=\alpha |k\rangle$ for some $k\in\{0,1,\cdots ,l-1\}$.

Let us consider the case $-sq^{l-1}\ne
\pm {\rm i}$. For $\alpha\not=0$ we get contradiction with the
commutation relations (1)--(3) applying them to the vector
$|l-1\rangle$. For $\alpha=0$ and $l\ge 2$ we get the one-dimensional
representation on the invariant subspace ${\bf C}|l-1\rangle$. For
$\alpha=0$ and $l=1$ we must move attention to the
vectors $|0\rangle, |-1\rangle,\cdots $  and the rest of the proof
is fulfilled by repeating the above and below arguments except that we
work with vectors with negative indices.

Now let us consider the case $-s q^{l-1}=\varepsilon {\rm i}$,
$\varepsilon=\pm 1$. Since $sq^{l-1}=-s^{-1}q^{-l+1}$, the equalities (34)
and (35) have no sense in this case.

For $\alpha \ne 0$
we get from equation $[l-1+j]_{q,s}=[l-1-j]_{q,s}$ (valid for all $j
\in {\bf Z}$) that $k=l-2$. By (27), (30) and the relation
$sq^{l-1}=-s^{-1}q^{-l+1}$ we have
$$
|l\rangle=R({\rm i} T_1 -s^{-1} q^{-l+3/2}
T_2)|l-1\rangle=\alpha |k\rangle=\alpha R({\rm i} T_1+s
q^{l-1/2} T_2)|l-1\rangle=0.
$$
This contradict the equality $|l\rangle =\alpha |k\rangle $, $\alpha \ne 0$.
For $\alpha=0$ we have from (1) and from $R({\rm i}
T_1-\varepsilon {\rm i} q^{1/2} T_2)|l-1\rangle=0$ that
$$
R(I)R(T_2) |l-1\rangle=\varepsilon \frac{q^2+1}{q^2-1} R(T_2)
|l-1\rangle={\rm i} [l]_{q,s} R(T_2) |l-1\rangle
$$
and we can redefine $|l\rangle:=R(T_2) |l-1\rangle$. From
(3) and (30) we have $0=R({\rm i} T_1+s q^{l+1/2} T_2) |l\rangle$.

If $|l\rangle$ is linear dependent on $|0\rangle, |1\rangle , \cdots ,
|l-1\rangle$, there exists $\beta\in {\bf C}$ such that
$|l\rangle=\beta |l-2\rangle$. As above,
for $\beta\ne 0$ we get the contradiction $\beta |l-1\rangle=0$ and
for $\beta=0$ we get a one-dimensional
representation since ${\bf C}|l-1\rangle$ is an invariant subspace.

If $|l\rangle$ is linear independent on the vectors $|0\rangle,
|1\rangle ,\cdots , |l-1\rangle$, we recursively
redefine $|j+1\rangle:=R({\rm i} T_1
s^{-1} q^{-j+1/2} T_2) |j\rangle$, $j=l,l+1,\cdots $. Then we
again consider 2 cases:

If there exists $\gamma \in {\bf C}$ such that $|l+l'\rangle=\gamma
|l-2-l'\rangle$ for some $l' \in \{1,2,\cdots ,l-2\}$ then we get
either one-dimensional representation on the invariant subspace
${\bf C}|l+l'-1\rangle$ (when $\gamma=0$) or a contradiction applying
(1)--(3) to $|l-2-l'\rangle$.

If $|l+l'\rangle$ is linear independent on $|0\rangle,|1\rangle ,\cdots ,
|l+l'-1\rangle$, then the representation is reducible since ${\rm lin}\,
\{|j\rangle,|j+1\rangle,\cdots \}$ is invariant subspace for any
$j$.

(b) $C \ne 0$. Consider first the case when $|l\rangle$ is
linear dependent on the linear independent vectors
$|0\rangle,|1\rangle ,\cdots , |l-1\rangle$. It means that
$|l\rangle=\alpha |k\rangle$ for some $k \in \{0,1,\cdots ,l-1\}$ and
for $\alpha \in {\bf C}$. If $\alpha=0$ we get a contradiction since
$$
0=|l\rangle=R({\rm i} T_1-s^{-1} q^{-l+3/2}
T_2)|l-1\rangle=R({\rm i} T_1+s q^{l+1/2}
T_2)|l\rangle=-Cq|l-1\rangle,
$$
implies $C=0$. For $\alpha\not=0$ we get a contradiction
applying (1)--(3) to the vector $|l-1\rangle$.

Now consider the case when the vectors $|j\rangle$, $j=0,1,2,\cdots $, are
linear independent. If there exists $m \in {\bf N}$ such that
the vector $|-m\rangle$ is linear
dependent on the linear independent vectors $|j\rangle$, $j \in
\{-m+1,-m+2,...\}$, we write $|-m\rangle=\beta |p\rangle$ for some
$\beta \in {\bf C}$ and $p>-m$. For $\beta=0$ we get a contradiction
similarly as in the above analogous cases. For $\beta =Cq$ and $p=-m+1$
we get the representation given (after some suitable rescaling of the basis)
by (22)--(26). For $\beta=-Cq$ and $p=-m+2$ we can
derive from (1) how operator $R(I)$ acts on the
linear independent vectors $|p\rangle$ and $R(T_2)|p-1\rangle$
and see that it cannot be diagonalized on this subspace.
Therefore, this case is unpossible.
For other values of $\beta$ and $p$ we get a contradiction applying
(1)--(3) to $|-m+1\rangle$.

Thus the only possible remaining case is when all the vectors
$|j\rangle$, $j \in {\bf Z}$ are linear independent. Using the
formulas (29), (32) and (33), in this case we
get (after some suitable rescaling of the basis) the representation
(18)--(20). Theorem is proved.
\medskip

It is clear from Theorem 2 that for $q\in {\bf R}$ the irreducible
$*$-representations of $U_q({\rm iso}_2)$ which
can be separated from representations of Theorem 1, are equivalent to
the irreducible $*$-representations from [7]. However, it is not seen
directly from formulas for representations
since operators of represenations in [7] are given
with respect to another basis than our one. Namely, the authors of [7]
diagonalize the operator $R(T_2)$ which corresponds to the shifts
in the group $ISO(2)$.
\medskip

\noindent {\bf ACKNOWLEDGEMENTS}
\medskip

We thank each others institutions for hospitality during mutual
visits. The research of M. H. and S. P. was supported by research
grants from GA Czech Republic. The research of A. U. K. was supported in
part by CRDF Grant UP1-309.
\medskip

\noindent {\bf REFERENCES}
\medskip

\noindent
1. L. L. Vaksman and L. I. Korogodskii, {\it Soviet Math. Dokl.} {\bf 39},
173 (1989).

\noindent
2. N. Ja. Vilenkin and A. U. Klimyk, {\it Representations of Lie Groups and
Special Functions}, vol. 1, Kluwer, Dordrecht, 1991.

\noindent
3. A. U. Klimyk, Preprint ITP--90--37E, Kiev, 1990.

\noindent
4. M. Odesski, {\it Funct. Anal. Appl.} {\bf 20}, No. 2, 78 (1986).

\noindent
5. Yu. S. Samoilenko and L. B. Turovska, in book {\it Quantum Groups and
Quantum Spaces}, Banach Center Publications, vol. 40, Warsaw, 1997,
pp. 21--40.

\noindent
6. M. Havl\'{i}\v cek, A. Klimyk and S. Po\v sta,
{\it J. Math.  Phys.}, {\bf 40} (1999), to be published.

\noindent
7. S. D. Silvestrov and L. D. Turowska, {\it J. Funct. Anal.}, {\bf 160},
79 (1998).

\noindent
8. A. M. Gavrilik and N. Z. Iorgov, Proc. of Second Intern. Conf.
{\it Symmetry in Nonlinear Mathematical Physics}, Kiev, 1997, 384.

\noindent
9. G. M. Bergman, {\it Adv. Math.} {\bf 29}, 178 (1978).

\noindent
10. A. Klimyk and K. Schm\"udgen, {\it Quantum Groups and Their
Representations}, Springer, Berlin, 1997.

\noindent
11. J. Dixmier, {\it Alg\'ebras Enveloppantes}, Gauthier--Villars, 1974.

\end{document}